\documentclass[a4paper,10pt]{amsart}
\usepackage{latexsym} 
\usepackage{amsfonts} 
\usepackage{amssymb} 
\usepackage{color}
\usepackage{marginnote}
\usepackage{amscd,enumerate}
\usepackage{epsfig}
\usepackage{amsmath}
\usepackage[all]{xy}
\font\pri=eufm10 at 11pt
\def\ppr#1{\hbox{\pri#1}}
\newcommand{\Z}{{\mathbb{Z}}}

\newtheorem{theo}{Theorem}

\newtheorem{co}{Corollary}

\newtheorem{re}{Remark}

\def\a{\alpha}

\def\fg{\mathfrak{L}}
\def\ft{\mathfrak{t}}

\def\ad{\textrm{ad}\;}
\title[Inner derivations of exceptional Lie algebras in prime characteristic]{Inner derivations of exceptional Lie algebras in prime characteristic}

\author[P. Alberca]{Pablo Alberca Bjerregaard}
\address{P. Alberca Bjerregaard: Departamento de Matem\'atica Aplicada, Escuela T\'ecnica Superior de Ingenieros Industriales, Universidad de M\'alaga. 29071 M\'alaga. Spain.}
\email{pgalberca@uma.es}
\author[D. Mart\'{\i}n]{Dolores Mart\'{\i}n Barquero}
\address{D. Mart\'{\i}n Barquero: Departamento de Matem\'atica Aplicada, Escuela T\'ecnica Superior de Ingenieros Industriales, Universidad de M\'alaga. 29071 M\'alaga. Spain.}
\email{dmartin@uma.es}

\author[C. Mart\'{\i}n]{C\'andido Mart\'{\i}n Gonz\'alez}
\address{C. Mart\'{\i}n Gonz\'alez:  Departamento de \'Algebra Geometr\'{\i}a y Topolog\'{\i}a, Fa\-cultad de Ciencias, Universidad de M\'alaga, Campus de Teatinos s/n. 29071 M\'alaga. Spain.}
\email{candido@apncs.cie.uma.es}

\date{1 abril 2013}

\begin{document}

\maketitle

\begin{abstract}
It is well-known that  every derivation of a semisimple Lie algebra $\fg$ over an algebraically closed field $F$ with characteristic zero  is inner. 
The aim of this paper is to show what happens if the characteristic of $F$ is prime with $\fg$ an exceptional Lie algebra. 
We prove that if $L$ is a Chevalley Lie algebra of type $\{\ppr{g}_2,\ppr{f}_4,\ppr{e}
_6,\ppr{e}_7,\ppr{e}_8\}$ over a field of characteristic $p$ then 
the derivations of $L$ are inner except in the cases $\ppr{g}_2$ with $p=2$,
$\ppr{e}_6$ with $p=3$ and $\ppr{e}_7$ with $p=2$.
\end{abstract}

\section{Introduction}
This paper deals with the derivation algebra of the Chevalley Lie algebras of exceptional type $\{\ppr{g}_2,\ppr{f}_4,\ppr{e}
_6,\ppr{e}_7,\ppr{e}_8\}$. As it is well-known Chevalley constructed Lie algebras over arbitrary fields $F$ starting from 
any semisimple (finite-dimensional) Lie algebra over an algebraically closed field of characteristic zero. The key point was
to realize that on such algebras one can find a suitable basis whose structure constants are integers. Then, by an scalar extension
process he constructed those Lie algebras which nowadays are called Chevalley algebras. 

By gathering results of Zassenhaus (1939), Seligman (1979), Springer and Steinberg (we will give full details in the next section) one can
get convinced that any derivation of a Chevalley $F$-algebra of any of the types $\{\ppr{g}_2,\ppr{f}_4,\ppr{e}_6,\ppr{e}_7\}$ is inner 
for $\mbox{\rm char}(F)\geq 5$. Also the derivations of a Chevalley algebra of type $\ppr{e}_8$ with $\mbox{\rm char}(F)\geq 7$, are inner.

In \cite{Elf4}, the case $\ppr{f}_4$ over fields of characteristic $\ne 2$ is
considered. There, it is proved that all derivations are inner. 
In this paper we study the cases not covered previously:
the derivations of the algebras of type $\{\ppr{g}_2,\ppr{f}_4,\ppr{e}_6,\ppr{e}_7\}$ over fields of characteristic
$2$ or $3$, and algebras of type $\ppr{e}_8$ over fields of characteristic $2$, $3$ or $5$ (we overlap with \cite{Elf4} only in the case
$\ppr{f}_4$ for characteristic $3$ though our methodology in this case is different).

Returning to the reference \cite{Elf4}, the key idea in this paper is the formula stating that any derivation is the
sum of an inner derivation plus a derivation annihilating a fixed Cartan subgalgebra (\cite[Proof of Proposition 8.1]{Elf4}). 
We have proved that a certain extension of this idea is also true for Chevalley algebras of exceptional type: any derivation is the sum of an inner derivation
plus a derivation mapping a fixed Cartan subalgebra into the center. Thus, when the center is null, we recover the formula in \cite{Elf4}.
This formula reduces the computation of the dimension of the derivation algebra to the computation of the dimension of the vector
space $V$ of derivations mapping a fixed Cartan subalgebra into the center.

For the exceptional algebras of low rank $\ppr{g}_2$ and $\ppr{f}_4$, most of the computational algebra packages with support for Lie algebras gives
us directly the dimension of the derivation algebra. However 
this computational approach does not work for the algebras $\ppr{e}_i$, $i=6,7,8$ which have a higher rank.

\section{Previous definitions and results}

%%%%%%%%%%%%%%%%%%%%%%%%%%%%%%%%%%%%%%%COMIENZO DE EOM1 %%%%%%%%%%%%%%%%%%%%%%%%%
Let $k$ be an algebraically closed field $k$ of characteristic $0$.
Recall that the Cartan matrix of a finite-dimensional semi-simple Lie algebra $\fg$ over $k$  is a matrix
$$A= \Bigg( 2\frac{(\a_i,\a_j)}{(\a_j,\a_j)}\Bigg)_{i,j = 1,\dots,r}$$
where $\a_1,\dots,\a_r$ is some system of simple roots of $\fg$ with respect to a fixed
Cartan subalgebra $\ft$ and $(\;,\;)$ is the scalar product on the dual space of $\ft$ defined by the
Killing form on $\fg$. The entries $a_{ij} = 2 \frac{(\a_i,\a_j)}{(\a_j,\a_j)}$ of a Cartan matrix have the following properties:

\begin{equation}  \left.\begin{aligned}   a_{ii}=2;\quad a_{ij}\le 0 \ &\textrm{and}\ a_{ij} \in \Z\textrm{ for }i\ne j,\\
        a_{ij}=0\ &\;\;\Rightarrow\  a_{ji}=0.       \end{aligned}
  \right\}\qquad\tag{1}  
  \end{equation}
%%%%%%%%%%%%%%%%%%%%%%%%%%%%%%%%%%%%%%%Fin de EOM1 %%%%%%%%%%%%%%%%%%%%%%%%%%%%%%

The Cartan matrix is a key tool in the description of $\fg$ by generators and relations. Denote by $\Phi$ the root system of $\fg$ relative to $K$. As it is known, there exist in $\fg$ a basis
$\{v_\a\}_{\alpha\in\Phi}\cup\{h_i\}_{i=1}^r$ whose structure constants are in $\Z$:
\begin{enumerate}\label{inadd}
\item $[h_i,h_j]=0$, for any $i,j\in\{1,\dots,r\}$.
\item~\label{miraque}  $[h_i,v_\a]=\langle\a,\a_i\rangle v_\a\in\Z v_\a$ for $\a\in\Phi, 1\le i\le r$ and $\langle\a,\beta\rangle:=2\frac{(\a,\beta)}{(\beta,\beta)}$.
\item  $[v_\a,v_{-a}]=h_\a$ a $\Z$-linear combination of  $h_1,\ldots,h_r$.
\item  If $\a$ and $\beta$ are independent roots, $\beta-r\a,\ldots,\beta+q\a$ the $\a$-string through $\beta$ then $[v_\a,v_\beta]=0$ if $q=0$
while $[v_\a,v_\beta]=\pm(r+1)v_{\a+\beta}$ if $\a+\beta\in\Phi$.
\end{enumerate}
Such a basis $\{v_\a\}_{\alpha\in\Phi}\cup\{h_i\}_{i=1}^r$ is called a {\em Chevalley basis} of $\fg$ and its existence is proved for instance in
\cite[Proposition, p. 146 and Theorem, p. 147]{Hum2}. 
%\marginpar{Falta n\'{u}mero del teorema}
 
\def\F{\mathbf F}
Following \cite[p. 149]{Hum2}, given a Chevalley basis $\{v_\a\}_{\alpha\in\Phi}\cup\{h_i\}_{i=1}^r$ of the semisimple Lie algebra $\fg$ over $k$ we can consider its linear envelope $\fg_\Z$
which is a Lie algebra over $\Z$ in the obvious sense. Now if $\F_p$ is the field of integers module the prime $p$ we can construct the Lie $\F_p$-algebra
$\fg_{\F_p}:=\fg_\Z\otimes_\Z\F_p$. The multiplication table of this algebra in a suitable basis is the one for $\fg$ reduced module $p$.
Now, if $F$ is any extension field of $\F_p$ we can construct the Lie algebra (over $F$) given by $\fg_F:=\fg_{\F_p}\otimes_{\F_p} F$. 
This algebra inherits both basis and Lie algebra structure from $\fg_{\F_p}$. In this way we obtain from each couple $(\fg,F)$ a Lie algebra $\fg_F$ over $F$,
called a Chevalley algebra,
which depends up to isomorphism only on $\fg$ and the field $F$.
If $\fg$ is simple of Dynkin diagram $D$,
the we will say that $\fg_F$ is a Lie algebra of type $D$. A common notation of this algebra will be $D(F)$ (for instance $\ppr{f}_4(\F_3)$ or $\ppr{e}_8(\F_2)$).

Let $\fg$ be a Lie algebra over a field $F$. If $x\in \fg$, we have $\mbox{ad}(x):\fg\to \fg$ with $\mbox{ad}(x)(y):=[x,y]$. It is easy to check that $\mbox{ad}(x)\in \mbox{Der}(\fg)$. We will denote $\mbox{ad}(\fg):=\{\mbox{ad}(x):x\in\fg\}\subset \mbox{Der}(\fg)$. In fact, 
\begin{equation}
\mbox{ad}(\fg)\triangleleft \mbox{Der}(\fg).
\end{equation}

\begin{theo}
Let $F$ be an algebraically closed field of characteristic zero. If $\fg$ is a finite-dimensional semisimple Lie algebra over $F$, then 
\begin{equation}
\mbox{\rm ad}(\fg)= \mbox{\rm Der}(\fg)
\end{equation} 
\end{theo}

If $\fg$ is finite-dimensional, we have the {\em Killing form} of $\fg$:
\begin{equation}
k:\fg\times\fg\to F,\ k(x,y):=\mbox{\rm tr}(\mbox{ad}(x),\mbox{ad}(y)),
\end{equation}
and we have the following results
%\begin{theo}[Cartan]
%Let $\fg$ be a finite-dimensional  Lie algebra over a field of characteristic zero. If the Killing form is non-degenerate then $\fg$ is semisimple.
%\end{theo}\smallskip

\begin{theo}[Zassenhaus \cite{Za}]~\label{ZA}
If $\fg$ is  a finite-dimensional  Lie algebra with non-degenerate Killing form, then every derivation is inner. That is to say 
\begin{equation}
 \mbox{\rm ad}(\fg)= \mbox{\rm Der}(\fg).
\end{equation}
\end{theo}

Let $\fg$ be a Chevalley Lie algebra of exceptional type $\{\ppr{g}_2,\ppr{f}_4,\ppr{e}_6,\ppr{e}_7,\ppr{e}_8\}$. A very fundamental result (\cite[p.49] {Hum1}) 
in what follows, due to Seligman and independently to Springer and Steinberg, says that (relative to a Chevalley basis) 
\begin{eqnarray}
\det{k}=2^\alpha 3^\beta,\ \mbox{ if } \fg \mbox{ is of type } \ppr{g}_2,\ppr{f}_4,\ppr{e}_6\mbox{ or  }\ppr{e}_7,~\label{detk1}\\
\det{k}=2^\alpha 3^\beta5^\gamma,\ \mbox{ if } \fg \mbox{ is of type } \ppr{e}_8~\label{detk2},
\end{eqnarray}
for some nonzero scalars $\alpha,\beta$ and $\gamma$. It is easy now to obtain
\smallskip

\begin{co}~\label{cor7}
Let $\fg$ be a Chevalley $F$-algebra of any of the types $\{\ppr{g}_2,\ppr{f}_4,\ppr{e}_6,\ppr{e}_7\}$, then if  $\mbox{\rm char}(F)\geq 5$ every derivation
is inner. If $\fg$ is of type  $\ppr{e}_8$ and $\mbox{\rm char}(F)\geq 7$, every derivation in $\fg$ is inner.
\end{co}

We have defined $\mbox{\rm ad}:\fg\to\mbox{Der}(\fg)$ as $x\mapsto \mbox{ad}(x)$.
An important fact is that   
\begin{equation}~\label{dimad}
\dim \mbox{ad}(\fg)=\dim \fg-\dim\ker\mbox{ad}=\dim \fg-\dim Z(\fg),
\end{equation}
where $Z(\fg)$ is the {\em center} of the Lie algebra. If, for example,  $Z(\fg)=0$ then $\mbox{ad}$ is a monomorphism and $\dim\fg=\dim\mbox{ad}(\fg)\leq \dim\mbox{Der}(\fg)$. Now if $\dim\fg=\dim\mbox{Der}(\fg)$  then $\mbox{ad}(\fg)=\mbox{Der}(\fg)$, and every derivation will be inner.\medskip

Take now a simple Lie algebra $\fg$ over an algebraically closed field $k$ of characteristic zero. Fix a Chevalley basis $\{v_\a\}\cup\{h_i\}$ of
$\fg$. The different
$h_i$ generate a Cartan subalgebra $H$ and each $v_\alpha$ the root space $\fg_\alpha$.
Consider now an arbitrary field $F$ and consider the Chevalley algebra $\fg_F$ introduced above. The properties of the new algebra $\fg_F$ may have changed dramatically, for instance $\fg_F$
may happen to be non-simple. Also while $\fg$ is simple, the algebra $\fg_F$ may have a nonzero center.  In any case 
$\fg_F=H_F\oplus(\oplus_\alpha (\fg_\alpha)_F)$ where $H_F$ is the $F$-linear span of the $h_i$'s and 
 $(\fg_\alpha)_F$ is the one-dimensional space $Fv_\alpha$. 
\medskip

\begin{re}\label{nota}\rm
 If $\alpha\colon H\to k$ is a nonzero root of $\fg$ then we can define a root $\bar\a$ of $\fg_F$ taking into account that $H_F$ is generated as $F$-vector
 space by the $h_i$'s of the Chevalley basis $\{v_\a\}_{\alpha\in\Phi}\cup\{h_i\}_{i=1}^r$. We must realize that $\a(h_i)\in\Z$ by (\ref{miraque}) of the definition of Chevalley basis. Thus, we define $\bar\a\colon H_F\to F$ as the $F$-linear extension $\bar\a(\sum_i\lambda_i h_i):=\sum_i\lambda_i \a(h_i)$.
 Since $\a$ is nonzero $\bar\a$ is also a nonzero.
 Furthermore it is straightforward that for any $h\in H_F$ one has $[h,v_\a]=\bar\a(h)v_\a$ in $\fg_F$.
\end{re}

The result in the remark below is a corollary of the main result in \cite{Dieu}.

\begin{re}\label{notados}\rm
 Let $\fg_F$ be a Chevalley algebra of exceptional type over a field $F$ of prime characteristic. 
 Then the center of $\fg_F$ is zero except for $\ppr{e}_6$ in characteristic $3$ and $\ppr{e}_7$ in characteristic $2$.
 In this cases the center $Z$ is one-dimensional and $Z\subset H_F$.
\end{re}

The following results takes the idea of \cite[Proposition 8.1]{Elf4} to whom we are emdebted. Though in \cite{Elf4} the result appears in a particular
context, the idea can be extended to the following setting.
\medskip

\begin{theo}~\label{mtheo}  Let $\fg$ be a semisimple Lie algebra over an algebraically closed field $k$
of characteristic zero and let $\fg_F$ be the corresponding Chevalley $F$-algebra (where $F$ is a field of prime characteristic).
Denote by $V$ the $F$-vector space of all derivations $d$ of $\fg_F$ such that $d(H_F)\subset Z:=Z(\fg_F)$, then $\mbox{\rm Der}(\fg_F)=\mbox{\rm ad}(\fg_F)+V$ hence we have the following Hesse diagram of subspaces 
\[
\xymatrix{
& \mbox{\rm Der}(\fg_F) & \cr
  \mbox{\rm ad}(\fg_F) \ar[ur]   & & \ar[ul] V  \cr
  & \ar[ur]\ar[ul]\mbox{\rm ad}(\fg_F)\cap V  & 
 }
%\caption{Rhombus.}
\]
\end{theo}
Proof. Consider first the Cartan decomposition $\fg=H\oplus(\oplus_{\a}\fg_\a)$ of $\fg$ with relation to $H$. 
As it is known this decomposition is a (fine) grading on $\fg$
whose zero component is $H$ and the rest of the homogeneous components are the root spaces $\fg_\a$. 
Moreover, it is a group grading with grading group $\Z^r$ being $r=\dim H$.
Next we take the
induced $\Z^r$-grading on the Chevalley algebra (zero component $H_F$ and homogeneous components $(\fg_\a)_F$) which will play
an important role in the forthcoming argument. Now we consider the $\Z^r$-grading on $\mbox{\rm Der}(\fg_F)$ such that a derivation $d$ is of degree $z$
if and only if $d$ takes the component of degree $w$ of $\fg_F$ to the component of degree $w+z$ (for any $z,w\in\Z^r$).
In order to prove the formula $\mbox{\rm Der}(\fg_F)=\mbox{\rm ad}(\fg_F)+V$ it suffices to prove that for any homogeneous derivation $d$
we can decompose $d$ as $d=\ad(x)+v$ for some $x\in\fg_F$ and $v\in V$. We analyze first the case in which $d$ is of degree $0$. Then $d(H_F)\subset H_F$
and $d(v_\alpha)=\lambda_\a v_\a$ for some scalar $\lambda_\a\in F$. So starting from $[h,v_\a]=\a(h)v_\a$ and applying $d$ we get
$[d(h),v_\a]+\lambda_\a [h,v_\a]=\a(h)\lambda_\a v_a$ (for arbitrary $\a$). Thus $[d(h),v_\a]=0$ and $d(h)\in Z(\fg_F)$ so that in this case $d\in V$.
Now assume that $d$ is a derivation of degree $z\ne 0$ and that $d(H_F)\ne 0$ (if $d$ annihilates $H_F$ then $d\in V$).  If the component of degree $z$ of
$\fg_F$ is $(\fg_\alpha)_F$ then for any $h\in H_F$ we have $d(h)=\lambda(h)v_\a$ where $\lambda\colon H_F\to F$ is a linear map. 
If we now take arbitrary elements $h,k\in H_F$, since $[h,k]=0$ we get 
$$\lambda(h)[v_\a,k]+\lambda(k)[h,v_\a]=0,$$
implying $\lambda(h)\bar\a(k)=\lambda(k)\bar\a(h)$. The fact that $d$ does not annihilates $H_F$ implies $\lambda\ne 0$ hence 
there is some $h\in H_f$ such that $\bar\a=c\lambda$ being $c=\lambda(h)^{-1}\bar\a(h)$. Take into account also that $c\ne 0$ since $\bar\a\ne 0$
as pointed out in Remark \ref{nota}. Finally the reader can check that $d+c^{-1}\ad(v_\a)$ is a derivation of $\fg_F$ in $V$ (more precisely
it annihilates $H_F$). 

\begin{co}\label{couno}~\label{fdimD}
If $Z$ is the center of $\fg_F$ then 
 $$\dim \mbox{\rm Der}(\fg_F)=\dim_k \fg+\dim V-\dim_k H.$$
 (we write $\dim$ for $\dim_F$).
\end{co}
Proof. We know that $\dim \mbox{\rm ad}(\fg_F)=\dim\fg_F-\dim Z=\dim_k\fg-\dim Z$. 
On the other hand $\mbox{\rm ad}(\fg_F)\cap V=\{ad(x)\colon [x,H_F]\subset Z\}$ therefore 
$\dim \mbox{\rm ad}(\fg_F)\cap V=\dim \{x\colon [x,H_F]\subset Z\}-\dim Z$.
Next we prove that $\{x\colon [x,H_F]\subset Z\}=H$. Take $x=h+\sum_\a\lambda_\a v_a$ with $h\in H_F$ and $\lambda_\a\in F$,
satisfying $[x,H_F]\subset Z$.
Then $\sum_\a\lambda_\a\a(k)v_\a\in Z$ for any $k\in H$.
Since $Z\subset H_F$ we conclude $\sum_\a\lambda_\a\a(k)v_\a=0$ hence if some $\lambda_\a\ne 0$ then $\a(H)=0$ and so $\a=0$
a contradiction.

So $\dim \mbox{\rm Der}(\fg_F)=\dim_k\fg-\dim Z+\dim V-\dim_k H+\dim Z$ whence
the result.\medskip

\begin{theo}\label{referee} {\rm (Referee's private communication)}
Under the hypothesis of Theorem \ref{mtheo}, if the root system of $\fg$ is simply laced, then $\dim(V)= \dim(H)$
and so $$\dim(\ad(\fg_F))=\dim(\text{Der}(\fg_F))-\dim Z.$$
Thus if $\fg_F$ is centerless,
all its derivations are inner.
\end{theo}
Proof.
We consider a Chevalley basis $\{v_\alpha\}_{\alpha\in\Phi}\cup\{h_i\}_{i=1}^r$ of $\fg$ as described in section 2.
We write $[v_\alpha , v_\beta] = n_{\alpha,\beta} v_{\alpha+\beta}$ where $n_{\a,\beta}\in F$ .
Since the root system is simply laced (for example of type $E$),
we have $n_{\alpha,\beta} = \pm 1$, and in particular it is nonzero modulo $p$.
Let $d$ be a derivation of $\fg_F$ such that $d(H_F) \subset Z$. Then,
the root spaces of $\fg_F$ with respect to $H_F$ are
$d$-invariant. Now in type $E$, the root spaces are either $1$-dimensional, or
$2$-dimensional. The latter case only happens in characteristic $2$, where $\alpha$
and $-\alpha$ agree. In any case, there are scalars $\lambda_{\alpha,i}$ and $\mu_{\alpha,i}$, such that 
$$\begin{cases}
   d(v_\a)=\lambda_{\a,1}v_\a+\lambda_{\a,2}v_{-\a}\cr
   d(v_{-\a})=\mu_{\a,1}v_\a+\mu_{\a,2}v_{-\a}.
  \end{cases}$$
Let now $\a = \a_i$ be a simple root, then by applying $d$ to $[v_\a , v_{-\a} ] = h_i$ we
get $d(h_i) = (\lambda_{\a,1} + \mu_{\a,2} )h_i$ . But also $d(h_i) \in Z$. Observe that no $h_i$ lies
in the centre (it would force a row of the Cartan matrix to be zero, for
example). It follows that $d(h_i ) = 0$ and $\mu_{\a,2} = -\lambda_{\a,1}$. In particular, we
get $d(H_F ) = 0$.

Let $\beta$ be another simple root such that $\a + \beta \in \Phi$. Applying $d$ to
the equality $[v_\a , v_\beta ] = n_{\a,\beta} v_{\a+\beta}$ we get $d(v_{\a+\beta}) = (\lambda_{\a,1} + \lambda_{\beta,1} )v_{\a+\beta}$ . Also applying
$d$ to $[v_{-\a}, v_{\a+\beta}] = n_{-\a,\a+\beta} v_\beta$ we get $d(v_\beta ) = \lambda_{\beta,1} v_\beta$, that is, $\lambda_{\beta,2} = 0$. Similarly,
we get $\mu_{\beta,1} = 0$. It follows that $d(v_{\a} ) = \lambda_{\a} v_{\a}$ , $d(v_{\a} ) = -\lambda_{\a} v_{-\a}$ for all
simple roots $\a$.

Conversely, fix a scalar $\lambda_\a$ for all simple roots $\a$, and define a linear map
$d\colon \fg_F \to \fg_F$ by $d(H_F) = 0$ and $d(v_\gamma ) = av_\gamma$, where $a = \sum_\a c_\a \lambda_\a$, being
$\gamma = \sum_\a c_\a \a$. Then by using the multiplication table of $\fg_F$ we see that $d$
is a derivation of $\fg_F$ . The conclusion is that the space $V$ of all derivations of $\fg_F$ mapping $H_F$ into $Z$ 
has dimension exactly $r$ (the rank of the root
system).

Thus, the formula in Corollay \ref{fdimD} reduces to 
 $\dim \mbox{\rm Der}(\fg_F)=\dim_k \fg$ and since $\dim(\fg)-\dim(Z)=\dim(\ad(\fg_F))$ we finally get
 $\dim(\ad(\fg_F))=\dim \mbox{\rm Der}(\fg_F)-\dim(Z)$. In particular when $Z=0$ all the derivations are inner.
 
\begin{co}
  For $i\in\{6,7,8\}$ the derivations of Chevalley algebras of exceptional type $\ppr{e}_i$  are all inner except in the cases $\ppr{e}_6$ in characteristic $3$ 
  and $\ppr{e}_7$ in characteristic $2$.
\end{co}
Proof. Let $\fg_F$ be a Chevalley algebra of type $\ppr{e}_i$ with $i=6,7$ or $8$. Taking into account Remark \ref{notados} we have $Z(\fg_F)=0$ except for
$\fg_F$ of type $\ppr{e}_6$ and the characteristic of the ground field is $3$, or $\fg_F$ of type $\ppr{e}_7$ and the characteristic of the ground field is $2$.
Then applying the formula in Theorem \ref{referee} we get the announced result.

\section{Exceptional lie algebras $\ppr{g}_2$ and $\ppr{f}_4$}
%\marginnote{\tt Esta secci\'on esta por re-escribir eliminado computaci\'on}

With these two algebras, $\ppr{g}_2$ and $\ppr{f}_4$, we have used the software GAP to compute the 
dimension of their derivation algebras. All the following routines have run reasonably   rapidly in a personal computer.

\subsection{Lie algebra  $\ppr{g}_2$}

Let us consider a Chevalley Lie algebra $\ppr{g}_2$. If the characteristic of the base field is other than 2 or 3 then every derivation is inner (Corollary \ref{cor7}). We present a computational approach to the problem and extend this result to the cases not covered. We are going to  prove the following result:

\begin{theo}
Let $\fg=\ppr{g}_2$ over a field $F$ of prime characteristic  $p$. Every derivation in $\ppr{g}_2$ is inner if, and only if, $p\not=2$. Thus
\begin{equation}
\mbox{\rm ad}(\ppr{g}_2)=\mbox{\rm Der}(\ppr{g}_2)\mbox{ if }\mbox{\rm char}(F)\not= 2.
\end{equation}
If the characteristic is $2$ then $\dim\mbox{\rm Der}(\ppr{g}_2)=21$ while $\dim\ppr{g}_2=14$.
\end{theo}

\textsc{Proof.} If the characteristic of the base field is $p\geq 5$ we can use Corollary \ref{cor7} and (\ref{detk1}), and then every derivation is inner.
 Let us now consider $p=3$, since $\dim\mbox{\rm Der}(\ppr{g}_2(F))=\dim\mbox{\rm Der}(\ppr{g}_2(\mathbf{F}_3))$, we may take without loss of generality $F=\mathbf{F}_3$. We use the software GAP. The next lines 

%%%%%%%%%%%%%%%%
\def\peq{\fontsize{9}{9}\selectfont}

\bigskip

{\peq 
\hrule
\begin{center}
\begin{minipage}{10cm}
gap$>$ F:={\bf GF}(3);\\
gap$>$ L:={\bf SimpleLieAlgebra}("G",2,F);\\
$<${\em Lie algebra of dimension 14 over GF(3)}$>$
\end{minipage}
\end{center}
\hrule
}
\bigskip

\noindent define the base field $F$ of characteristic $3$ and the Lie algebra $\ppr{g}_2$ over $F$, which is 14-dimensional. We can also compute de center of the Lie algebra by doing

\bigskip

{\peq 
\hrule
\begin{center}
\begin{minipage}{10cm}
gap$>$ {\bf LieCenter}(L);\\
$<${\em Lie algebra of dimension 0 over GF(3)}$>$
\end{minipage}
\end{center}
\hrule
}

\bigskip

\noindent which is 0-dimensional. We determine the dimension of the Lie algebra $\mbox{\rm Der}(\ppr{g}_2)$:

\bigskip

{\peq  
\hrule
\begin{center}
\begin{minipage}{10cm}
gap$>$ B:={\bf Basis}(L);\\
{\em CanonicalBasis( $<$Lie algebra of dimension 14 over GF(3)$>$ )}\\
gap$>$ {\bf Derivations}(B);\\
$<${\em Lie algebra of dimension 14 over GF(3)$>$}
\end{minipage}
\end{center}
\hrule
}

\bigskip

\noindent and then $\dim\mbox{\rm Der}(\ppr{g}_2)=14$. We have 
$
\dim\mbox{\rm ad}(\ppr{g}_2)=\dim\ppr{g}_2-\dim Z(\ppr{g}_2)=14-0=14=\dim\mbox{\rm Der}(\ppr{g}_2)$.
Thus every derivation is inner also if $\mbox{\rm char}(F)=3$ as  we have confirmed that 
\begin{equation}
 \mbox{\rm ad}(\ppr{g}_2)=\mbox{\rm Der}(\ppr{g}_2), \mbox{\rm char}(F)=3,
\end{equation}
in this case by using a dimensional reasoning.

Now, if $p=2$, we do first
\bigskip

{\peq 
\hrule
\begin{center}
\begin{minipage}{10cm}
gap$>$ F:={\bf GF}(2);\\
gap$>$ L:={\bf SimpleLieAlgebra}("G",2,F);\\
$<${\em Lie algebra of dimension 14 over GF(2)}$>$
\end{minipage}
\end{center}
\hrule
}

\bigskip

\noindent in order to work with $\ppr{g}_2$ with $\mbox{\rm char}(F)=2$, which is also 14-dimensional.  We have now to define a basis and then we can compute the dimension of   $\mbox{\rm Der}(\ppr{g}_2)$:
\medskip

{\peq 
\hrule
\begin{center}
\begin{minipage}{10cm}
gap$>$ B:={\bf Basis}(L);\\
{\em CanonicalBasis( $<$Lie algebra of dimension 14 over GF(2)$>$ )}\\
gap$>$ {\bf Derivations}(B);\\
$<${\em Lie algebra of dimension 21 over GF(2)$>$}
\end{minipage}
\end{center}
\hrule
}

\bigskip

\noindent and it has dimension 21. Then   
\begin{equation}
\mbox{\rm ad}(\ppr{g}_2)\subsetneq \mbox{\rm Der}(\ppr{g}_2),\ \mbox{\rm char}(F)=2,
\end{equation}
because   $\dim \mbox{\rm ad}(\ppr{g}_2)=\dim \ppr{g}_2=14\neq 21$, where we have used that the center is again 0-dimensional, which can be confirmed   with the next code line: 

%\vfill\eject

\bigskip

{\peq 
\hrule
\begin{center}
\begin{minipage}{10cm}
gap$>$ {\bf LieCenter}(L);\\
$<${\em Lie algebra of dimension 0 over GF(2)}$>$
\end{minipage}
\end{center}
\hrule
}
\bigskip
 
Finally, we present a table that summarizes the previous computations:

\bigskip
{\footnotesize
\begin{center}
\begin{tabular}{||c|c|c|c|c|c|c||}
\hline
$\fg$ & char.& $\dim \fg$ & $\dim Z(\fg)$ & $\dim \mbox{\rm ad}(\fg)$ & $\dim \mbox{\rm Der}(\fg)$ & $\mbox{\rm ad}(\fg)=\mbox{\rm Der}(\fg)$ \\
\hline
$\ppr{g}_2$ & 2 & 14 & 0 & 14 & 21 & no \\
\hline
$\ppr{g}_2$ & 3 & 14 & 0 & 14 & 14 & yes  \\
\hline
\end{tabular}
\end{center}
}
\bigskip

\hfill $\square$

\subsection{Lie algebra $\ppr{f}_4$}

Let us work with the Chevalley Lie algebra $\ppr{f}_4$. If the characteristic of the base field is other than 2 or 3, we can use 
Corollary \ref{cor7} to conclude that every derivation is inner. In spite of the fact that A. Elduque and M. Kochetov have recently proved that the result 
is also true if the characteristic of the base field is 3 (see \cite[Proposition 8.1]{Elf4}) we present an extension that confirms 
this result is also true at any prime characteristic.

If the characteristic of the base field is 2 or 3 we have repeated the previous strategy using GAP with the results summarized in the following table:

\bigskip
{\footnotesize
\begin{center}
\begin{tabular}{||c|c|c|c|c|c|c||}
\hline
$\fg$ & char.& $\dim \fg$ & $\dim Z(\fg)$ & $\dim \mbox{\rm ad}(\fg)$ & $\dim \mbox{\rm Der}(\fg)$ & $\mbox{\rm ad}(\fg)=\mbox{\rm Der}(\fg)$ \\
\hline
$\ppr{f}_4$ & 2 & 52 & 0 & 52 & 52 & yes \\
\hline
$\ppr{f}_4$ & 3 & 52 & 0 & 52 & 52 & yes  \\
\hline
\end{tabular}
\end{center}
}
\bigskip
\noindent   Then we have the next result:

\begin{theo}
If $\fg=\ppr{f}_4$ over a field $F$ of prime characteristic then 
\begin{equation}
\mbox{\rm ad}(\ppr{f}_4)=\mbox{\rm Der}(\ppr{f}_4),
\end{equation}
that is to say, every derivation of $\ppr{f}_4$ is inner. 
\end{theo}

\section{Further considerations}

We summarize all the results and the software involved at the following table:

%\vskip 1cm

\begin{table}[h]
{\footnotesize
\begin{center}
\begin{tabular}{||c|c|c|c|c|c|c||}
\hline
$\fg$ & char.& $\dim$ & $Z(\fg)$ & $ \mbox{\rm ad}(\fg)$ & $  \mbox{\rm Der}(\fg)$ & inner \\
\hline
$\ppr{g}_2$ & 2 & 14 & 0 & 14 & 21 & no   \\
\hline
$\ppr{g}_2$ & 3 & 14 & 0 & 14 & 14 & yes   \\
\hline
$\ppr{f}_4$ & 2 & 52 & 0 & 52 & 52 & yes  \\
\hline
$\ppr{f}_4$ & 3 & 52 & 0 & 52 & 52 & yes  \\
\hline
$\ppr{e}_6$ & 2 & 78 & 0 & 78 & 78 & yes  \\
\hline
$\ppr{e}_6$ & 3 & 78 & 1 & 77 & 78 & no  \\
\hline
$\ppr{e}_7$ & 2 & 133 & 1 & 132 &  133 & no \\
\hline
$\ppr{e}_7$ & 3 & 133 & 0 & 133 & 133 & yes  \\
\hline
$\ppr{e}_8$ & 2 & 248 & 0 & 248 & 248 & yes  \\
\hline
$\ppr{e}_8$ & 3 & 248 & 0 & 248 & 248 & yes  \\
\hline
$\ppr{e}_8$ & 5 & 248 & 0 & 248 & 248 & yes  \\
\hline
\end{tabular}
\end{center}
}
\caption{Dimensions and results.}
\end{table}

And as a corollary	

\begin{theo}
 The derivations of Chevalley algebras of exceptional type are all inner except in the cases $\ppr{g}_2$ in characteristic $2$,
 $\ppr{e}_6$ in characteristic $3$ and
 $\ppr{e}_7$ in characteristic $2$.
\end{theo}

\section*{Acknowledgments}
%%%%%%%%%%%%%%%%%%%%%%%%%%%%%%
%%
All the authors have been partially supported by the Spanish MEC and Fondos FEDER through project MTM2010-15223,  by the Junta de Andaluc\'{\i}a and Fondos FEDER, jointly, through projects FQM-336, FQM-02467 and FQM-3737.

The author thankfully acknowledges the computer resources, technical
expertise and assistance provided by the SCBI (Supercomputing and
Bioinformatics) center of the University of M\'{a}laga.

\end{document}